\documentclass{amsart}
\usepackage[cp1251]{inputenc}
\usepackage[english]{babel}
\usepackage{amsmath}
\usepackage{amssymb}
\usepackage{amsfonts}
\newtheorem{thm}{Theorem}
\begin{document}
\thispagestyle{empty}

 \title[VECTOR TRANSFORM OPERATORS ]{  VECTOR TRANSFORM OPERATORS FOR PIECE-WISE HARMONIC FUNCTIONS
 }%

\author{O.Yaremko, Y. Parfenova}

\address{Oleg Yaremko, Yulia Parfenova
\newline\hphantom{iii}Penza State University,
\newline\hphantom{iii}str. Lermontov, 37, 
\newline\hphantom{iii} 440038, Penza, Russia}
\email{yaremki@mail.ru}

\maketitle {\small

\begin{quote}
\noindent{\bf Abstract. } The vector transform operators are investigated; these operators are used at the solution of boundary
value problems in piecewise homogeneous spherically symmetric areas. In particular,
examples of transformation operators for vector boundary value problems are given
for third vector boundary value problem in the unit circle and for the Dirichlet problem in the unit circle.
\medskip


\end{quote} }

\section{Introduction.}

 The transform operator is the operator which translates the decision of one
problem of mathematical physics in the decision of other problem. Transform
operators of K. Weierstrass, S. Poisson, N. Ja. Sonin are known .

Our interest are the special cases of transform operators concerning
different boundary value problems for the same equation. By the way of
illustration it is considered the first and third boundary value problems of
Dirichlet for the Laplase equation in unit circle:
\[
\Delta \tilde {u}=0,
\]
\[
u|_s =f\left( \varphi \right),
\]
\[
h \tilde{u}+ {\tilde {u}|}_{r=1} =f\left( \varphi \right)
\quad
h>0.
\]

It is possible to prove that the transform operator $P:\tilde {u}\to u$, as
shown I.I.Bavrin in work [1], has form

\[
u\left( x \right)=\int\limits_0^1 {\varepsilon ^{h-1}\tilde {u}\left(
{\varepsilon x} \right)d\varepsilon }.
\]

In the present work transform operators are constructed in a vector case
which formally turns out replacement of function $u$ by a vector function
$u$, and number $h$ by a matrix $H$.

The case of transform operators is studied, the transform operators
connecting the decision of the first boundary value problem with internal
conditions of interface on sphere and the decision of a problem of Dirichlet
are constructed.

Statement of the first regional problem with interface conditions:

\[
\Delta u_k =0,\,x\in V_k ;\;k=1,...,n+1.
\]

The edge conditions are

\[
\Gamma _0 \left[ {u_1 } \right]=f_0 \left( \eta \right),
\quad
\eta \in S_0,
\]
Where
\[
u_k =\left( {\begin{array}{l}
 u_{k1} \\
 u_{k2} \\
 \vdots \\
 u_{km} \\
 \end{array}} \right),
\quad
f_0 =\left( {\begin{array}{l}
 f_{01} \\
 f_{02} \\
 \vdots \\
 f_{0m} \\
 \end{array}} \right).
\]

There are the non-uniform contact on the hypersurfaces conjugation

\[
S_k , S_k =\left\{ {\eta =\left( {\eta _1 ,...,\eta _N } \right):\left\|
\eta \right\|=r_k } \right\}:
\]
\[
\Gamma _{j1}^k \left[ {u_k } \right]-\Gamma _{j2}^k \left[ {u_{k+1} }
\right]=f_{jk} \left( \eta \right);
\quad
\eta \in S_k ;
\quad
k=1,...,n,
\quad
j=1,2,
\]

Where
\[
f_{jk} =\left( {\begin{array}{l}
 f_{jk1} \\
 f_{jk2} \\
 \vdots \\
 f_{jkm} \\
 \end{array}} \right),
\]

Here

\[
\Delta u_k =\frac{1}{r^{N-1}}\frac{\partial }{\partial r}\left(
{r^{N-1}\frac{\partial u_k }{\partial r}} \right)+\frac{1}{r^2}\Delta _\eta
u_k ,
\]
$\Delta _\eta $
 is Laplase operator on the sphere $S_0 $; $\Gamma _0 ,\,\Gamma
_{j1}^k ,\,\Gamma _{j2}^k \left( {j=1,2;\,k=1,...,n} \right)$ are posed
operators, permutable with operator $L_0 =\sum\limits_{i=1}^N {x_i
\frac{d}{dx_i }} $.

The transform operator looks like:
\[
u=P_0 \left[ {\tilde {u}_0 } \right]+\sum\limits_{j=1}^2
{\sum\limits_{s=1}^n {P_{js} \left[ {\tilde {u}_{js} } \right]} }
\]
where $P_0 ,\,P_{jq} $ - vector transform operators.

Kelvin reflection method for solutions of  mathematical physics boundary value problem with symmetric boundary from served
as the basis for method of operators progresses in mathematical physics, complex analysis, harmonic analysis[2].
 In this paper operator method  is developed for vector problems of heterogeneous pattern mathematical physics.

\section{The common boundary value problem for the Laplase equation in unit sphere with non-uniform internal conditions of interface.}

Let  $B_{n} $ be piesewise homogeneous unit ball of $R^{N} $ :

     $$B_{n}^{} =\mathop{\bigcup }\limits_{i=1}^{n+1} V_{i} ;V_{i} =\left\{\; x\in R^{N} :r_{i} <\left\| x\right\| <r_{i-1} \right\}\; ;i=1,...,n+1,$$

 $$B_{n} =S_{0} \times I_{n}^{+} ,  S_{0} =\left\{\eta \in R^{N} :\left\| \eta \right\| ^{2} =1\right\},$$

$$I_{n}^{+} =\left\{r:r\in \bigcup \limits _{j=1}^{n+1}\left(r_{j} ,r_{j-1} \right);0<r_{0} \le 1,r_{n+1} =0,r_{j+1} <r_{j} ,j=0,...,n \right\}.$$

     Let us considere a problem about construction set separate Laplace combined equations solution, bounded on $B_{n}^{} $
      \begin{equation}
      \displaystyle \Delta u_{k} =0,x\in V_{k} ;k=1,...,n+1;
      \end{equation}
by boundary conditions
\begin{equation}
\displaystyle \Gamma _{0} \left[u_{1} \right]=f_{0} \left(\eta \right),\quad \eta \in S_{0}
\end{equation}
Where
\[
u_k =\left( {\begin{array}{l}
 u_{k1} \\
 u_{k2} \\
 \vdots \\
 u_{km} \\
 \end{array}} \right),
\quad
f_0 =\left( {\begin{array}{l}
 f_{01} \\
 f_{02} \\
 \vdots \\
 f_{0m} \\
 \end{array}} \right).
\]
There are the heterogeneous contact on conjunction hypersurfaces conjugations

 $S_{k} $ , $S_{k} =\left\{\; \eta =\left(\eta _{1} ,...,\eta _{N} \right):\left\| \eta \right\| =r_{k}^{} \; \right\}$ :
\begin{equation}
      \displaystyle \Gamma _{j1}^{k} \left[u_{k} \right]-\Gamma _{j2}^{k} \left[u_{k+1} \right]=f_{jk}^{} \left(\eta \right)\; \quad ;\eta \in S_{k} ;k=1,...,n,\; \, j=1,2,
 \end{equation}
where
\[
f_{jk} =\left( {\begin{array}{l}
 f_{jk1} \\
 f_{jk2} \\
 \vdots \\
 f_{jkm} \\
 \end{array}} \right),
\]
Here

$$\displaystyle \Delta u_{k} =\frac{1}{r^{N-1} } \frac{\partial }{\partial r} \left(r^{N-1} \frac{\partial u_{k} }{\partial r} \right)+\frac{1}{r^{2} } \Delta _{\eta } u_{k}, $$
 $\Delta _{\eta } $ is the Laplase operator on sphere  $S_{0} $ ; $\Gamma _{0} ,\Gamma _{j1}^{k} ,\Gamma _{j2}^{k}       \; \left(j=1,2;\; \, k=1,...,n\right)$  are permutable with $L_0 =\sum\limits_{i=1}^N {x_i \frac{d}{dx_i }} $  defined operators.

  \section{Method of influence function.}
       Fourier transform on sphere $S_{0} $  with nonseparated  variables
[3] reduces the problem (1)-(3) to form: find the separate combined differential equations solution

      $$\frac{1}{r^{N-1} } \frac{d}{dr} (r^{N-1} \frac{d\tilde{u}_{k,l} }{dr} )-l\left(l+N-2\right)\frac{1}{r^{2} } \tilde{u}_{k,l} =0\; ;l\in Z,r_{i} <r<r_{i-1} $$
by boundary conditions

      $$\left. \Gamma _{0} \left[\tilde{u}_{1,l} \right]\right|_{r=r_{0} } =r_{0}^{N-1} \tilde{f}_{0,l} $$
There are the heterogeneous contact conditions in points  joint $r=r_{k} $

      $$\displaystyle \Gamma _{j1,l}^{k} \left[\tilde{u}_{k,l} \right]-\Gamma _{j2,l}^{k} \left[\tilde{u}_{k+1,l} \right]=r_{k}^{N-1} \tilde{f}_{jk,l} ,k=1,...,n\; ;j=1,2.$$

Assign formula by immediate checking:

$$\tilde{u}_{j,l} \left(r\right)=H_{j,1,l}^{*} \left(r,r_{0} \right)\; \left(\begin{array}{c} {0} \\ {1} \end{array}\right)r_{0}^{N-1} \tilde{f}_{0,l} +\sum \limits _{s=1}^{n}H \, _{j,s,l}^{*} \left(r,r_{s} \right)\; r_{s}^{N-1} \left(\begin{array}{c} {\tilde{f}_{1s,l} } \\ {\tilde{f}_{2s,l} } \end{array}\right),$$
Here  $H_{k,s,l}^{*} =H_{k,s,l}^{*} \left(r,\rho \right)$  are matrix-valued $m\times m$ functions, that defined by formulas:

When  {\it k }{\it \textless }{\it  s}

$$H_{k,s,l}^{*} =\left(\varphi _{k,l} \left(r\right)\, \, -\psi _{k,l} \left(r\right)\mathop{\psi _{1,l}^{-1} }\limits^{0\quad }  \mathop{\varphi _{1,l}^{} }\limits^{0\quad } \right) \left(1\quad 0\right)\; \, \Omega _{s,l}^{-1} \left(\rho \right),r_{k} <r<r_{k-1} ,\; \, r_{s} <\rho <r_{s-1} ,$$

When ${ k \textgreater  s }$

$$H_{k,s,l}^{*} =-\psi _{k,l} \left(r\right)\left(\mathop{\psi _{1,l}^{-1} }\limits^{0\quad }  \left(\begin{array}{cc} {\mathop{\varphi _{1,l}^{} }\limits^{0} } & {\mathop{\psi _{1,l}^{} }\limits^{0} } \end{array}\right)\; \, \Omega _{s,l}^{-1} \left(\rho \right)\right),r_{k} <r<r_{k-1} ,\; \, r_{s} <\rho <r_{s-1} ,$$

When  {{\it k} = {\it s}}

      $$H_{k,s,l}^{*} =\left\{\begin{array}{l} {\left(\varphi _{k,l} \left(r\right)\, \, -\psi _{k,l} \left(r\right)\mathop{\psi _{1,l}^{-1} }\limits^{0\quad }  \mathop{\varphi _{1,l}^{} }\limits^{0\quad } \right) \left(1\quad 0\right)\; \, \Omega _{s,l}^{-1} \left(\rho \right),} \\ {r_{k-1} <r<\rho <r_{k} ,} \\ {-\psi _{k,l} \left(r\right)\left(\mathop{\psi _{1,l}^{-1} }\limits^{0\quad }  \left(\begin{array}{cc} {\mathop{\varphi _{1,l}^{} }\limits^{0\quad } } & {\mathop{\psi _{1,l}^{} }\limits^{0\quad } } \end{array}\right)\; \, \Omega _{s,l}^{-1} \left(\rho \right)\right),} \\ {r_{k-1} <\rho <r<r_{k} \; .} \end{array}\right. $$

Let us define matrix-valued functions  $$\varphi _{n+1,l} \left(r\right)=r^{l}\cdot E ;
\psi _{n+1,l} \left(r\right)=r^{-\left(l+N-2\right)}\cdot E ;l\in Z,$$
where $E$ is identity matrix.
 Order  {\it n} pairs of functions  $\left(\varphi _{k,l} ,\psi _{k,l} \right)\; ,k=1,...,n$  are founded from the recurrence equations.

$$\Gamma _{j1}^{k} \left(\varphi _{k,l} ,\psi _{k,l} \right)=\Gamma _{j2}^{k} \left(\varphi _{k+1,l} ,\psi _{k+1,l} \right)\; ,k=1,...,n;j=1,2.$$

Let us use the following notations

$$\left. \mathop{\varphi }\limits^{0} _{1,l} =\Gamma _{0} \left[\varphi _{1,l} \left(r\right)\right]\; \, \right|_{\; r=r_{0} } ,\left. \mathop{\psi }\limits^{0} _{1,l} =\Gamma _{0} \left[\psi _{1,l} \left(r\right)\right]\; \, \right|_{\; r=r_{0} } ;l\in Z,$$

$$\Gamma _{ij}^{k} \left(\varphi _{k,l} \left(r\right),\psi _{k,l} \left(r\right)\right)=\left(\varphi _{ij,l}^{k} \left(r\right),\psi ^{k} _{ij,l} \left(r\right)\right)\; ;i,j=1,2\; ;k=1,...,n\; ;l\in Z,$$

      $$\Omega _{k,l} \left(\rho \right)=\left(\begin{array}{cc} {\varphi _{11,l}^{k} \left(\rho \right)} & {\psi _{11,l}^{k} \left(\rho \right)} \\ {\varphi _{12,l}^{k} \left(\rho \right)} & {\psi _{12,l}^{k} \left(\rho \right)} \end{array}\right)\; .$$

Matrix-valued functions $H_{k,s,l}^{*} \left(r,\rho \right)$  correctly defined,  if following conditions satisfied:

i)  when $l\to \infty $,  matrix sequences $\alpha _{0,l} ,\alpha _{j1,l}^{k} ,\alpha _{j2,l}^{k} $, that defined by formulas $\Gamma _{0} \left[r^{l} \right]=\alpha _{0,l}  r^{l} ,\Gamma _{ij}^{k} \left[r^{l} \right]=\alpha _{ij,l}^{k}  r^{l} ;i,j=1,2;k=1,...,n;l\in Z$ ,

have growth which no more than power-mode, and
$$\det M_{kj,l} =\det\left(\begin{array}{cc} {\alpha _{1j,l}^{k} } & {\alpha _{1j,-l}^{k} } \\ {\alpha _{2j,l}^{k} } & {\alpha _{2j,-l}^{k} } \end{array}\right)\ne 0;k=1,...,n;j=1,2;l\in Z,$$

ii) for every {\it l} $\in Z$  following inequality are valid $\det \Omega _{k,l} \left(\rho \right)\ne 0\; ;$ $k=1,...,n$ ; $\mathop{\psi _{1,l} }\limits^{0} \ne 0.\; $

    Receive formula for problem solutions(1)-(3) by returning fo Fourier original:

$$u_{j} \left(r\xi \right)=\frac{1}{\omega _{N} } \int _{S_{0} }^{} \left(H_{j,1}^{} \left(r,r_{0} ,\left\langle \eta ,\xi \right\rangle \right)\; \left(\begin{array}{c} {0} \\ {1} \end{array}\right)r_{0}^{N-1} f_{0} \left(\eta \right)+\right. $$
\begin{equation}
      \displaystyle \left. +\sum \limits _{s=1}^{n}H \, _{j,s}^{} \left(r,r_{s} ,\left\langle \eta ,\xi \right\rangle \right)\; r_{{\rm s}}^{N-1} \left(\begin{array}{c} {f_{1s} \left(\eta \right)} \\ {f_{2s} \left(\eta \right)} \end{array}\right)\right)dS_{0},
  \end{equation}
Where

$$H\, _{j,s}^{} \left(r,r_{s} ,\left\langle \eta ,\xi \right\rangle \right)=\sum _{l=0}^{\infty }\frac{2l+N-1}{N-1} C_{l}^{(N-1)/2} \left(\left\langle \eta ,\xi \right\rangle \right)H_{j,s,l}^{*} \left(r,r_{s} \right)\;  $$
-  $\omega _{N} $  - ({\it N}-1)-  $l\to \infty $ -dimensional volume of unit sphere $S_{0} $  from  $R^{N} $ ;   $C_{l}^{(N-1)/2} $ - are Gegenbauer polynomials [1], ${\left\langle {\eta ,\xi } \right\rangle }$ - is scalar product of vectors ${\eta ,\xi }$

    \section{Transform operators}
     Let  $\hat{u}_{0} ,\hat{u}_{jk} $ are harmonic vector-functions in the unit ball $B_{0}^{} $, $B_{0}^{} =\left\{\; x\in R^{N} :\left\| x\right\| ^{2} <1\right\}\; $ and continuous on  $\overline{B}_{0}^{} $ Boundary values of that functions are vectors $f_{0} \left(\eta \right),f_{jk} \left(\eta \right)$ respectively. Let us define vector transform operators  $P_{0} $ ,  $P_{jq} $  by using rules:

if
$$u_{0} \left(r\xi \right)=\sum \limits _{k=1}^{n+1}\chi \left(V_{k} \right) u_{0k} \left(r\xi \right),$$

$$\displaystyle u_{0k}^{} \left(r\xi \right)=\int _{S_{0} }^{}H_{k,1}^{} \left(r,r_{0} ,\left\langle \eta ,\xi \right\rangle \right)\left(\begin{array}{c} {0} \\ {1} \end{array}\right)\; r_{0}^{N-1} \;  \hat{u}_{0} \left(r\eta \right) dS_{0}, $$

then  $P_{0} :\hat{u}_{0} \to u_{0} $ ;

similarly, if
 $$u_{jq} \left(r\xi \right)=\sum \limits _{k=1}^{n+1}\chi \left(V_{k} \right)\, \; u_{jq,\; k} \left(r\xi \right)\,  $$

$$\displaystyle u_{jq,k}^{} \left(r\xi \right)=\int _{S_{0} }^{}H_{k,q,l} \left(r,r_{q} ,\left\langle \eta ,\xi \right\rangle \right)\left(\begin{array}{c} {\delta _{1j} } \\ {\delta _{2j} } \end{array}\right)\; r_{q}^{N-1}  \; \hat{u}_{jq} \left(r\eta \right) dS_{0}, $$

then  $P_{jq} :\hat{u}_{jq} \to u_{jq} ,j=1,2;q=1,...,n$ ,
where  $\chi \left(V_{k} \right)\, $ is a characteristic function in $V_{k} $ $\displaystyle \delta _{ij} $ is Kronecker symbol.
\[
\chi \left( {V_k } \right)=\left\{ {\begin{array}{l}
 1,\;r\xi \in V_k , \\
 0,\;r\xi \notin V_k \\
 \end{array}} \right.,
\quad
\delta _{ij} =\left\{ {\begin{array}{l}
 1,\;i=j, \\
 0,\;i\ne j. \\
 \end{array}} \right.
\]

      Basic formula (4) can be brought to the following form by using $P_{0} $ ,  $P_{jq} $  vector transform operators:
$$\displaystyle u=P_{0} \left[\hat{u}_{0} \right]+\sum _{j=1}^{2}\sum _{s=1}^{n}P_{js}   \left[\hat{u}_{js} \right].$$
     \begin{thm}If the conditions of existence i)-ii) vector-valued function $H_{k,s,l}^{*} \left(r,\rho \right),$ satisfied, then transformation operator $P_{0} $   $\left(P_{jq} \right)$  compare $\hat{u}_{0} $   harmonic vector function  $\left(\hat{u}_{jq} \right),$ in  the $B_{0}$   homogeneous ball and piesewise harmonic vector function  $u_{0} $   $\left(u_{jq} \right)$ in  piesewise homogeneous ball   $B_{n} ,$ The components $u_{0,k} $   $\left(u_{jq,k} \right)$  of that function continued in $V_{k} ,$ ball layer  and satisfy the boundary condition (2) and conjugation condition (3).
      \end{thm}

   \emph{   Example }1.      Transform operator  $P_{0} $  for the third vector boundary value problem  $\left. Hu+E\frac{\partial u}{\partial n} \right|_{S_{0} } =\left. \hat{u}\right|_{s_{0} } $ in unit cercle have form
$$P_{0} :\hat{u}\to u ,  u\left(x\right)=\int _{0}^{1}\varepsilon ^{H-E}  \hat{u}\left(\varepsilon x\right)d\varepsilon. $$
      Here
\[
\varepsilon ^{H-E}=e^{(H-E)\ln \varepsilon },
\]
 $E$ is identity matrix, $H$ is symmetric and positive-definite matrix [4].

     \emph{ Example }2. Transform operator  $P_{0} $  for Dirichlet problem {\it  }in the unit cercle with the internal conjunction conditions

         $$u^{-} \left(\eta \right)=u^{+} \left(\eta \right),\quad \quad K\frac{\partial }{\partial n} u^{-} \left(\eta \right)=\frac{\partial }{\partial n} u^{+} \left(\eta \right),\; \left|\eta \right|=r$$
have form  $P_{0} :\hat{u}\to u$ ,
\[
u\left( x \right)=\left\{ {\begin{array}{l}
 \sum\limits_{j=0}^\infty {\begin{array}{l}
 \left[ {\left( {E-K} \right)\left( {E+K} \right)^{-1}} \right]^j\left(
{\hat {u}\left( {xr^{2j}} \right)-\left( {E-K} \right)\left( {E+K}
\right)^{-1}\hat {u}\left( {\frac{x}{\left| x \right|^2}r^{2j+2}} \right)}
\right),\; \\
 \quad \quad \quad \quad \quad \quad \quad \quad \quad \quad \quad \quad
\quad \quad \quad \quad \quad \quad \quad \quad \quad \quad \quad \quad
r<\left| x \right|<1, \\
 \end{array}} \\
 2K\left( {E+K} \right)^{-1}\sum\limits_{j=0}^\infty {\left[ {\left( {E-K}
\right)\left( {E+K} \right)^{-1}} \right]^j\hat {u}\left( {xr^{2j}}
\right),\;} \left| x \right|<r. \\
 \end{array}} \right.
\]
Where $\left( {E+K} \right)^{-1}$ is inverse matrix, $u^{-} \left(\eta \right),u^{+} \left(\eta \right)$ are limit of function $u=u\left(x\right)$ when $x\to h$ from without and from within respectively.
Similarly, $$\frac{\partial }{\partial n} u^{-} \left(\eta \right), \frac{\partial }{\partial n} u^{+} \left(\eta \right)\quad are \quad limits \quad{u}'_r \left( x \right),$$ when $x\to h$ from without and from within respectively.
  \subsection{Transform operators in half-plane}
Method of transformation operators is used to solve the problem [1].
Necessary definitions from [3], [4], [5]. The direct $J:\hat {f}\to f$ f
and inverse $J^{-1}:f\to \hat {f}$ transformation operators are set
equalities:
\[
f(x)=\int\limits_{-\infty }^\infty {\varphi (x,\lambda )\left(
{\int\limits_{-\infty }^\infty {e^{-i\lambda \xi }\hat {f}(\xi )d\xi } }
\right)d\lambda ,}
\]
\[
\hat {f}(x)=\int\limits_{-\infty }^\infty {e^{-i\lambda \xi }\left(
{\int\limits_{-\infty }^\infty {\varphi ^\ast (\xi ,\lambda )f(\xi )d\xi } }
\right)d\lambda .}
\]
Here $\varphi (x,\lambda ),\varphi ^\ast (x,\lambda )$ - are the
eigenfunctions [13], [14] of the direct and coupling Sturm--Liouville
problems for the Fourier operator in piecewise-homogeneous axis In.
Eigenfunction
\[
\varphi (x,\lambda )=\sum\limits_{k=2}^n {\theta (x-l_{k-1} )} \theta (l_k
-x)\varphi _k (x,\lambda )+
\]
\[
+\theta (l_1 -x)\varphi _1 (x,\lambda )+\theta (x-l_n )\varphi _{n+1}
(x,\lambda )
\]
is a solution of the system of separate differential equations
\[
\left( {a_m^2 \frac{d^2}{dx^2}+\lambda ^2} \right)\varphi _m (x,\lambda
)=0,x\in \left( {l_m ,l_{m+1} } \right);m=1,...,n+1,
\]
on the coupling conditions
\[
\left[ {\alpha _{m1}^k \frac{d}{dx}+\beta _{m1}^k } \right]\varphi _k
=\left[ {\alpha _{m2}^k \frac{d}{dx}+\beta _{m2}^k } \right]\varphi _{k+1}
,
\]
on the boundary conditions
\[
\left. {\varphi _1 } \right|_{x=-\infty } =0,\left. {\varphi _{n+1} }
\right|_{x=\infty } =0
\]
Similarly eigenfunction
\[
\begin{array}{l}
 \varphi ^\ast (\xi ,\lambda )=\sum\limits_{k=2}^n {\theta (\xi -l_{k-1}
)\theta (l_k -\theta )} \varphi ^\ast (\xi ,\lambda )+ \\
 \theta (l_1 -\xi )\varphi _1^\ast (\xi ,\lambda )+\theta (\xi -l_n )\varphi
_{n+1}^\ast (\xi ,\lambda ) \\
 \end{array}
\]
is a solution of the system of separate differential equations
\[
\left( {a_m^2 \frac{d^2}{dx^2}+\lambda ^2} \right)\varphi _m^\ast (x,\lambda
)=0,x\in (l_m ,l_{m+1} );m=1,...,n+1,
\]
With the coupling conditions
\[
\frac{1}{\Delta _{1,k} }\left[ {\alpha _{m1}^k \frac{d}{dx}+\beta _{m1}^k }
\right]\varphi _k^\ast =\frac{1}{\Delta _{2,k} }\left[ {\alpha _{m2}^k
\frac{d}{dx}+\beta _{m2}^k } \right]\varphi _{k+1}^\ast ,x=l_k ,
\]
where
\[
\Delta _{i,k} =\det \left( {\begin{array}{l}
 \alpha _{1i}^k \beta _{1i}^k \\
 \alpha _{2i}^k \beta _{2i}^k \\
 \end{array}} \right)k=1,...,n;i,m=1,2,
\]
on the boundary conditions
\[
\left. {\varphi _1 } \right|_{x=-\infty } =0,\left. {\varphi _{n+1} }
\right|_{x=\infty } =0
\]
Let for some $\lambda $ of the considered boundary value problems have
nontrivial solutions $\varphi (x,\lambda ),\varphi ^\ast (x,\lambda )$, in
this case the number $\lambda $ is called the eigenvalue [13], [14],
corresponding solutions $\varphi (x,\lambda ),\varphi ^\ast (x,\lambda )$ -
is called the eigenfunctions of the direct and coupling Sturm--Liouville
problems, respectively. In the further we shall adhere to the following
normalization of eigenfunctions:
\[
\varphi _{n+1} (x,\lambda )=e^{ia_{n+1}^{-1} x\lambda }.\varphi _{n+1}^\ast
(x,\lambda )=e^{-ia_{n+1}^{-1} x\lambda }.
\]

\section{Conclusion.}

The vector transform operators are investigated; these operators are used at the solution of boundary
value problems in piecewise homogeneous spherically symmetric areas in the article.
Further it is supposed to extend results of work to a case of two and more
internal conditions of interface.

\end{document}